\documentclass [12pt,paper] {article}
\usepackage[cp866]{inputenc}
\usepackage[english,russian]{babel}
\usepackage{amssymb,amsmath,amsfonts,latexsym,amsthm}
\date{}
\newtheorem {theorem}{\bf Теорема}
\newtheorem {corol}{\bf Следствие}
\newtheorem {lemma}{\bf Лемма}

\newtheorem {prop}{\bf Утверждение}
\newcommand {\proofr }{{\par\medskip\noindent \bf Доказательство. }}

\newcommand {\eproof }{\hfill $\blacktriangle$ \\ \medskip}

\def\p0{\parindent 0pt}
\title{{\bf
\normalsize  И.~Ю.~Могильных, Ф.~И.~Соловьева\\
\large О максимальных компонентах одного класса совершенных
двоичных кодов}
\thanks{Работа  выполнена при  финансовой поддержке  гранта Российский Научный Фонд 14-11-00555.}}

\begin{document}

\maketitle

\begin{abstract}
В работе показано существование широкого  класса расширенных
совершенных двоичных кодов Соловьевой длины 16, содержащих
$ij$-компоненты максимальной мощности.
\end{abstract}

\section{Введение}

Компоненты совершенных двоичных кодов изучались большим числом
исследователей. Эти объекты имеют весьма сложную структуру,
используя их, удалось решить множество открытых проблем, стоящих в
теории совершенных кодов, см. обзоры \cite{S2008Switch,
Sobsor2013}. Множество $K$ называется $i$-{\it компонентой}, если
множество векторов, покрываемых $K$ и $K+e_i$ совпадают. Далее
будем рассматривать лишь {\it неразложимые} компоненты, то есть
те, которые нельзя разбить на компоненты меньшей мощности. Будем
называть $i$-компоненту {\it минимальной} ({\it максимальной}),
если она имеет наименьшую (наибольшую) возможную мощность.
Известно, что минимальная компонента в пространстве всех двоичных
векторов ${\textbf F}^n$, 
$n=2^r-1$, единственна с точностью до изоморфизма и может быть
представлена формулой $\{(x,|x|,x): x\in {\textbf F}^{(n-1)/2}\}$.
Максимальной мощности $i$-компоненты были впервые построены в 1988
г. в работе \cite{S88}. Оказалось, что эти компоненты занимают
половину множества кодовых слов совершенного кода. Существование
совершенных кодов с $i$-компонентами различных мощностей
установлено в 1995 г. в работе \cite{AS95}. В 2001 г. в статье
\cite{S01} было доказано существование максимальных неизоморфных
$i$-компонент, принадлежащих различным совершенным кодам для любой
допустимой длины $n=2^{m} - 1, m > 3$. П.\,Р.\,Ж.\,Остергард и
О.\,Поттонен, см. \cite{OP,OPP}, перечислили все возможные размеры
$i$-компонент совершенных кодов длины 15 и, в частности,
установили, что размер любой компоненты для любого совершенного
кода длины 15 всегда кратен размеру минимальной компоненты. Более
того, они перечислили все коды с указанием состава мощностей
$i$-компонент, см. \cite{Pot}. В  2012 г. В.\,Н.\,Потапов
\cite{P2010} описал возможные размеры компонент
 мощностей, близких к минимальной, для любой допустимой длины.

К.\,Т.\,Фелпс и М.\,ЛеВан \cite{PL99}, используя конструкцию
\cite{S81} автора настоящей статьи, в 1999 г. доказали, что
существуют совершенные коды длины 15, которые невозможно получить
из кода Хэмминга методом свитчинга. Этот класс кодов состоит из
двух неизоморфных совершенных кодов длины 15, каждый из которых по
любой координате разбивается на две максимальные компоненты.

 В
\cite{OPP} были перечислены все (многошаговые) свитчинговые классы
совершенных кодов длины 15, их оказалось девять, один из которых
имеет максимальную мощность, равную 5819, что представляет собой
подавляющую часть от общего числа 5983 всех неэквивалентных
совершенных кодов длины 15. Среди этих девяти классов обнаружено
четыре спорадических класса, каждый состоящий только из одного
совершенного кода.

\section{Обозначения и необходимые утверждения}

Далее в статье будем использовать транзитивные разбиения ${\textbf
F}^n$, $n=2^m,\,\, m\geq 3$ на максимально непараллельные
расширенные коды Хэмминга, предложенные Д.С.Кро\-товым в
\cite{Kr}. Для  $n=8$ такое разбиение  (разбиение под номером 8)
было приведено К.Т.Фелпсом в  \cite{P2000}, поэтому при $n=8$ ниже
будем ссылаться на него как на разбиение Фелпса.

Следуя терминологии и определениям \cite{Kr}, будем рассматривать
расширенный код Хэмминга длины $n=2^m,\,\, m\geq 3$ как
совокупность подмножеств $X$ в конечном поле ${\textbf
F}=GF(2^m)$, элементы которого,
 занумерованные в порядке
возрастания степеней примитивного элемента $\alpha$, есть кодовые
координаты.

Для  
$\alpha^k\in {\textbf F}$, $p\in \{0, 1\}$ определим код 
$H^p_{\alpha^k}$ как совокупность подмножеств $X$ в  поле
${\textbf F}$, удовлетворяющих проверочным соотношениям:

$$\sum_{x\in X} 1=p,$$
$$\sum_{x\in X}(x + \alpha^k)^3=0.$$
Также определим код $\overline{H}$ как циклический код с
порождающим многочленом, являющимся минимальным для $\alpha$:

$$\sum_{x\in X} 1=0,$$
$$\sum_{x\in X} x=0.$$

Легко видно, что $H^0_{\alpha}$ -- расширенный код Хэмминга, а
$H^1_{\alpha}$ -- его класс смежности. Всюду далее будем опускать
верхний индекс у нечетновесового кода $H^1_{\alpha}$ и будем его
обозначать $H_{\alpha}$, отвечающий ему код Хэмминга -- через
$\overline{H}_{\alpha}$. В \cite{Kr} доказано, что коды
$H_{\alpha}$, $\alpha\in {\textbf F}$ образуют транзитивное
разбиение множества нечетновесовых слов длины $n$, причем
$\overline{H}_{\alpha}$ и $\overline{H}_{\beta}$ для любых
различных $\alpha$ и $\beta$ из ${\textbf F}$ имеют минимально
возможное пересечение, равное $2^{2^m-2m}$ (при $n$, равном 8,
пересечение состоит из 4-х слов).

Для этих разбиений будем изучать строение компонент в коде
Соловьевой из \cite{S81}: $$\cup_{l \in \{0,1,\ldots,n\}}C_l\times
C_{\pi(l)},$$ где $\pi$ --  произвольная подстановка на множестве
элементов $\{0,1, \ldots,n\}$.

\begin{lemma}\label{l1}
Для любого $k\in {\mathcal F}$ код $\overline{H}\cap
\overline{H}_{\alpha^k}$ есть код БЧХ с порождающим многочленом
$m_1(x)m_3(x)$.
\end{lemma}
\proofr В силу определения, для кодового слова $X$, принадлежащего
одновременно кодам $\overline{H}$ и $\overline{H}_{\alpha^k}$,
имеем следующие проверочные соотношения:

 $$ \sum_{x\in X}1=0, \,\,\,
  \sum_{x\in X}x=0, \,\,\,
  \sum_{x\in X}(x+\alpha^k)^3=0.$$

 Из последнего равенства вытекает

 $$\sum_{x\in X}(x^3+x^2\alpha^k+x\alpha^{2k}+\alpha^{3k})=$$
$$\sum_{x\in X} x^3+\alpha^k (\sum_{x\in X}x)^2+\alpha^{2k}\sum_{x\in X}x+\alpha^{3k}\sum_{x\in X}1,$$
что, с учетом первых двух соотношений, дает
$$\sum_{x\in X}x^3=0.$$
Другими словами, слово $X$ принадлежит коду БЧХ с порождающим
многочленом $m_1(x)m_3(x)$. Обратное вложение доказывается
аналогично. \eproof

\begin{lemma}\label{l2}
Пусть слово $X$ веса 4 принадлежит коду $\overline{H}_{\alpha^i}
\cap \overline{H}_{\alpha^j}$. Тогда $X_{\alpha^i}$ и
$X_{\alpha^j}$ не могут одновременно равняться 1.
\end{lemma}

\proofr Пусть
$X=e_{\alpha^i}+e_{\alpha^j}+e_{\alpha^k}+e_{\alpha^l}$. В силу
определения кодов $\overline{H}_{\alpha^i}$,
$\overline{H}_{\alpha^j}$, имеем следующие соотношения:
$$(\alpha^j+\alpha^i)^3+(\alpha^r+\alpha^i)^3+(\alpha^s+\alpha^i)^3=0,$$
$$(\alpha^j+\alpha^i)^3+(\alpha^r+\alpha^j)^3+(\alpha^s+\alpha^j)^3=0,$$
 сложив которые, получим
$$
(\alpha^r+\alpha^i)^3+(\alpha^s+\alpha^i)^3+(\alpha^r+\alpha^j)^3+(\alpha^s+\alpha^j)^3=0.
$$

Преобразовывая левую сторону последнего равенства, выводим
равенство

$$(\alpha^i+\alpha^j)(\alpha^r+\alpha^s)(\alpha^i+\alpha^j+\alpha^r+\alpha^s)=0,$$
откуда по Лемме \ref{l1} следует, что $X$ принадлежит коду БЧХ
\,\, $\overline{H}\cap \overline{H}_{\alpha^i}$ и, следовательно,
не может иметь вес 4. \eproof
   \begin{corol}\label{inter} \,
Пусть $n=8$. \, Для произвольных $i$,  $j$ справедливо
$\overline{H}_{\alpha^i}\cap \overline{H}_{\alpha^j}=$
$=\{0^8,1^8,x^0,x^1\}$, где $x^0$, $x^1$ -- векторы веса 4,
единичные координатные позиции которых разбивают множество
координат на 2 множества $\Pi_{i,j}^0$ и $\Pi_{i,j}^1$, причем
$\alpha^i$ и $\alpha^j$ принадлежат различным множествам.
   \end{corol}

\section{Структура компонент разбиения Фелпса}
Для кодового слова $X$ расширенного совершенного кода $D$ через
$N_{ij}(D,X)$ обозначим множество кодовых слов, находящихся на
расстоянии $4$ от $X$ и различающихся с кодовым словом $X$ в $i$-й
и $j$-й координатах. Назовем эти кодовые слова $ij$-{\it смежными}
с кодовым словом $X$. Совокупность кодовых слов расширенного
совершенного кода $D$, которые можно соединить путем с кодовым
словом $X$, каждая последовательная пара которых состоит из
$ij$-смежных кодовых слов, называется $ij$-{\it компонентой} кода
$D$  и обозначается $R_{ij}(D,X)$.

Вначале выясним структуру \, $ij$-смежных \, кодовых слов в коде
\, Соловьевой $\cup_{l \in \{0,\ldots,n\}}C_l\times C_{\pi(l)}$
для произвольной подстановки $\pi$. Рассмотрим $ij$-смежные c
$(X,Y)\in C_k\times C_{\pi(k)}$ кодовые слова кода $\cup_{l \in
\{0,\ldots,n\}}C_l\times C_{\pi(l)}$. Прежде всего это слова,
получающиеся из $ij$-смежных слов с $X$ в коде $C_k$: $\{(Z,Y):Z
\in N_{ij}(C_k,X)\}$. Затем рассмотрим
 код $C_l$, где $X+e_i+e_j\in C_l$. Несложно видеть, что $Y$, как и всякое другое слово
 на расстоянии 2 от $C_{\pi(l)}$ находится на расстоянии 2 ровно от $n/2$
 кодовых слов кода $C_l$, которые имеют вид $Y +e_{s_t}+e_{s'_t}$,
 $t=1,\ldots, n/2$, где $\cup_{t=1,\ldots,n/2} \{s_t,s'_t\}$
 образует разбиение координат $\{1,\ldots,n\}$. Отсюда получаем,
 что $N_{ij}(D_P,(X,Y))=N_{ij}(C_k,X)\times Y \cup \{(X+e_i+e_j,Y+e_{s_t}+e_{s'_t}):t \in \{1,\ldots,n/2\}\}$.
 Таким образом доказали следующее утверждение, ввиду которого
приходим к выводу, что исследование $ij$-компонент кода Соловьевой
напрямую
 связано с изучением структуры графа расстояний Хэмминга 2,
 индуцированного кодовыми словами различных двух кодов разбиения
 $P$:

\begin{prop}
Пусть $P=\cup_{l\in\{0,\ldots,n\}} C_{l}$ есть разбиение множества
нечетновесовых векторов длины $n$ на расширенные совершенные коды,
  $(X,Y)$ -- кодовое слово кода Соловьевой
$D_P=\cup_{l\in\{0,\ldots,n\}} C_{l}\times C_{\pi(l)}$, $X \in
C_k$. Тогда, eсли $X+e_i+e_j\in C_l$, то $N_{ij}(D_P, (X,Y
))=N_{ij}(C_k,X)\times Y  \cup \{(X+e_i+e_j,Z): Z\in C_{\pi(l)},
d(Z,Y)=2\}$.
\end{prop}

\begin{corol}\label{corolcom}
Пусть $P=\cup_{j\in\{0,\ldots,n\}} C_{j}$ -- разбиение множества
нечетновесовых векторов длины $n$ на расширенные совершенные коды,
  $(X,Y)$ есть кодовое слово кода
$D_P=\cup_{j\in\{0,\ldots,n\}} C_{j}\times C_{\pi(j)}$, $X \in
C_k$. Тогда $R_{ij}(C_k,X)\times Y \subset R_{ij}(D_P,(X,Y))$.
\end{corol}

\begin{prop}\label{redfac}
Пусть $n=8$, \, $i$ и $j$ -- произвольные координаты, $i\neq j$.
Представителей классов смежностей $\overline{H}_{\alpha^i}$ по
$\overline{H}_{\alpha^i}\cap \overline{H}_{\alpha^j}$ можно
выбрать среди векторов веса не более 4 из $ij$-компоненты кода
$\overline{H}_{\alpha^i}$, 
содержащей нулевой вектор.
\end{prop}
\proofr Пусть
$e_{\alpha^i}+e_{\alpha^j}+e_{\alpha^{i_1}}+e_{\alpha^{i_2}}$  и
$e_{\alpha^i}+e_{\alpha^j}+e_{\alpha^{i_3}}+e_{\alpha^{i_4}}$
попадают в один класс смежности $\overline{H}_{\alpha^i}\cap
\overline{H}_{\alpha^j}+X$, где $X\in \overline{H}_{\alpha^i}$. В
силу того, что сумма этих двух векторов имеет нули на позициях
$\alpha^i$ и $\alpha^j$ и в силу антиподальности кода БЧХ с
порождающим многочленом $m_1(x)m_3(x)$, получаем, что
$\overline{H}_{\alpha^i}\cap \overline{H}_{\alpha^j}$ содержит
вектор веса 4 с единицами в позициях $\alpha^i$ и $\alpha^j$, что
противоречит Лемме \ref{l2}. \eproof

\begin{lemma}\label{l3}
1. Для всяких векторов $X\in {H}_{\alpha^i}$, $Y\in
{H}_{\alpha^j}$, класс смежности $X+Y+\overline{H}_{\alpha^i}\cap
\overline{H}_{\alpha^j}$ содержит единственного лидера, который
имеет вес 2.

2. Для различных классов смежности $H_{\alpha^i}+H_{\alpha^j}$ по
$\overline{H}_{\alpha^i}\cap \overline{H}_{\alpha^j}$ лидеры
различны и множество лидеров $X+Y+\overline{H}_{\alpha^i}\cap
\overline{H}_{\alpha^j}$ по всем $X\in \overline{H}_{\alpha^i}$,
$Y\in \overline{H}_{\alpha^j}$ есть $\{e_r+e_s:r \in
\Pi_{i,j}^0,s\in \Pi_{i,j}^1\}$.
\end{lemma}
\proofr Пусть носители векторов веса  4   кода \,
$\overline{H}_{\alpha^i}\cap \overline{H}_{\alpha^j}$ \, есть
$\Pi_{i,j}^0=$
$=\{\alpha^i,\alpha^{i_1},\alpha^{j_1},\alpha^{k_1}\}$
 и
 $\Pi_{i,j}^1=\{\alpha^j,\alpha^{i_2},\alpha^{j_2},\alpha^{k_2}\}$.

1. Векторы $X'=X+e_{\alpha^i}$, $Y'=Y+e_{\alpha^j}$ принадлежат
кодам $\overline{H}_{\alpha^i}$ и $\overline{H}_{\alpha^j}$
соответственно. Утверждение леммы очевидно выполнено в случае,
когда $X'$ или $Y'$
 принадлежат $\overline{H}_{\alpha^i}\cap
\overline{H}_{\alpha^j}$. По Утверждению \ref{redfac} без
ограничения общности имеем:
$X'=e_{\alpha^i}+e_{\alpha^j}+e_{\alpha^{i_1}}+e_{\alpha^{i_2}}$,
 тогда
 $Y'=e_{\alpha^i}+e_{\alpha^j}+e_{\alpha^{j_1}}+e_{\alpha^{j_2}}$
 или
 $Y'=e_{\alpha^i}+e_{\alpha^j}+e_{\alpha^{i_1}}+e_{\alpha^{j_2}}$.

 В первом случае в силу Утверждения \ref{inter} класс смежности $X'+Y'+\overline{H}_{\alpha^i}\cap
\overline{H}_{\alpha^j}=e_{\alpha^i}+e_{\alpha^j}+e_{\alpha^{j_1}}+e_{\alpha^{j_2}}+e_{\alpha^{i_1}}+e_{\alpha^{i_2}}+\overline{H}_{\alpha^i}\cap
\overline{H}_{\alpha^j}$ может иметь только один вектор веса 2, а
именно $e_{\alpha^{k_1}}+e_{\alpha^{k_2}}$. Во втором случае класс
смежности $X'+Y'+\overline{H}_{\alpha^i}\cap
\overline{H}_{\alpha^j}=e_{\alpha^i}+e_{\alpha^j}+e_{\alpha^{i_2}}+e_{\alpha^{j_2}}+\overline{H}_{\alpha^i}\cap
\overline{H}_{\alpha^j}$ содержит единственный вектор веса 2, а
именно $e_{\alpha^{i}}+e_{\alpha^{k_2}}$.

2. Очевидно, что различные пары классов смежности
$X+\overline{H}_{\alpha^i}\cap \overline{H}_{\alpha^j}$ и
$Y+\overline{H}_{\alpha^i}\cap \overline{H}_{\alpha^j}$ имеют
различных лидеров для класса $X+Y+\overline{H}_{\alpha^i}\cap
\overline{H}_{\alpha^j}$. В противном случае получаем противоречие
с тем, что  факторизуем по подпространству
$\overline{H}_{\alpha^i}\cap \overline{H}_{\alpha^j}$. Так как
классов смежности $H_{\alpha^i}+H_{\alpha^j}$ по
$\overline{H}_{\alpha^i}\cap \overline{H}_{\alpha^j}$ всего 16, а
их лидеры необходимо принадлежат множеству $\{e_r+e_s:r \in
\Pi_{i,j}^0,s\in \Pi_{i,j}^1\}$, получаем, что лидеры исчерпывают
это множество.

 \eproof

\begin{corol}\label{corolfactor} \,
Пусть \, $n=8$ и \, $k, \, r, \, s$ -- произвольные \, числа из
множества $\{-\infty,0,1,\ldots,6\}$, $r\ne s$, векторы $X$,
$X'\in {H}_{\alpha^k}$ находятся на расстоянии 4 и принадлежат
некоторой $(\alpha^r, \alpha^s)$-компоненте кода ${H}_{\alpha^k}$.
Тогда векторы $X+e_{\alpha^r}+e_{\alpha^s}$,
$X'+e_{\alpha^r}+e_{\alpha^s}$ не принадлежат одному и тому же
коду ${H}_{\alpha^l}$ для некоторого $l$.

\end{corol}
\proofr Пусть $X$ и $X'$ таковы, что
$X+e_{\alpha^r}+e_{\alpha^s}$, $X'+e_{\alpha^r}+e_{\alpha^s}\in
{H}_{\alpha^l}$  для некоторого $l$. Возможно два случая:
$X+X'=e_{\alpha^r}+e_{\alpha^s}+e_{\alpha^{i_1}}+e_{\alpha^{j_1}}$
и
$X+X'=e_{\alpha^{i_1}}+e_{\alpha^{j_1}}+e_{\alpha^{i_2}}+e_{\alpha^{j_2}}$.

Пусть
$X+X'=e_{\alpha^r}+e_{\alpha^s}+e_{\alpha^{i_1}}+e_{\alpha^{j_1}}$
и принадлежит $\overline{H}_{\alpha^l}\cap
\overline{H}_{\alpha^k}$. С другой стороны, вектор
$X+X'=e_{\alpha^r}+e_{\alpha^s}+e_{\alpha^{i_1}}+e_{\alpha^{j_1}}\notin
\overline{H}_{\alpha^l}\cap \overline{H}_{\alpha^k}$, так как по
Лемме \ref{l3} векторы $X$ и $X+e_{\alpha^r}+e_{\alpha^s}$ не
могут принадлежать кодам $H_{\alpha^k}$ и $H_{\alpha^l}$
соответственно (поскольку вектор $e_{\alpha^r}+e_{\alpha^s}$ не
принадлежит $\overline{H}_{\alpha^l}\cap
\overline{H}_{\alpha^k}$).

Второй случай доказывается аналогично. \eproof

Определим граф $G_{k,l}=(V,E)$ следующим образом: множество вершин
графа есть $V=H_{\alpha^k}\cup H_{\alpha^l}$, множество ребер есть
 $E=\{(X,X'):X\in H_{\alpha^k}, X'\in H_{\alpha^l}, d(X,X')=2\}$

\begin{lemma}\label{l4graph}
Пусть $n=8$ \, и \, $l, \, k$ -- произвольные различные числа из
множества $\{-\infty,0,1,\ldots,6\}$. Тогда граф $G_{k,l}$ связен.
\end{lemma}
\proofrНапомним, что в силу Следствия \ref{inter} восемь кодовых
координат можно разбить на два множества
$\Pi_{l,k}^0=\{l,i_1,j_1,k_1\}$ и $\Pi_{l,k}^1=\{k,i_2,j_2,k_2\}$.
Согласно Лемме \ref{l3} достаточно доказать, что всякая пара
вершин $X$,
$X'=X+e_{\alpha^l}+e_{\alpha^{i_1}}+e_{\alpha^{j_1}}+e_{\alpha^{k_1}}$
из одного класса смежности кода $H_{\alpha^k}$ по подкоду
$\overline{H}_{\alpha^k}\cap \overline{H}_{\alpha^l}$ соединена
путем в графе $G_{k,l}$. Так как рассматриваемые коды являются
расширенными совершенными, то для всякого кодового слова $X$
одного кода и любого $l$ найдется ровно одно кодовое слово в
другом коде, отличающееся от $X$ в $\alpha^{l}$-й позиции и
какой-то еще позиции, которую  обозначим $\alpha^{l*}$. В силу
Леммы \ref{l3} элементы $\alpha^l$ и $\alpha^{l*}$ не могут
одновременно принадлежать $\Pi_{l,k}^0$ или $\Pi_{l,k}^1$.
Построим путь от $X$ до $X'$ в графе $G_{k,l}$. Пусть $i'\in
\Pi_{l,k}^0\setminus l$:\\
 $X$, $X+e_{\alpha^l}+e_{\alpha^{l*}}$,
$X+e_{\alpha^l}+e_{\alpha^{l*}}+e_{\alpha^{i'}}+e_{\alpha^{i'*}}$,
$X+e_{\alpha^l}+e_{\alpha^{l*}}+e_{\alpha^{i'}}+e_{\alpha^{i'*}}+e_{\alpha^{l*}}+e_{\alpha^{(l*)*}}$,
$X+e_{\alpha^l}+e_{\alpha^{i'}}+e_{\alpha^{i'*}}+e_{\alpha^{(l*)*}}+e_{\alpha^{i'*}}+e_{\alpha^{(i'*)*}}=
X+e_{\alpha^l}+e_{\alpha^{i'}}+e_{\alpha^{(l*)*}}+e_{\alpha^{(i'*)*}}$.

По Лемме \ref{l3} элементы $\alpha^{(l*)*}, \alpha^{(i'*)*}\in
\Pi_{l,k}^0$ и отличаются от $\alpha^{l}$ и $\alpha^{i'}$ в силу
того, что рассматриваемые коды имеют кодовое расстояние 4. Отсюда
немедленно имеем, что
$X+e_{\alpha^l}+e_{\alpha^{i'}}+e_{\alpha^{(l*)*}}+e_{\alpha^{(i'*)*}}=X'$.
 \eproof

Пару координат кода $\cup_{l\in
\{-\infty,0,1,\ldots,6\}}H_{\alpha^l}\times H_{\alpha^{\pi(l)}}$
назовем {\it однородной}, если она  принадлежит первой или второй
половине кодовых координат одновременно.

\begin{theorem}
Пусть \, $P=\cup_{l\in \{-\infty,0,1,\ldots,6\}}H_{\alpha^l}$ \,
есть разбиение Фелпса.  \, Тогда код $D_P=\cup_{l\in
\{-\infty,0,1,\ldots,6\}}H_{\alpha^l}\times H_{\alpha^{\pi(l)}}$
имеет две максимальные компоненты по любой паре однородных
координат для произвольной перестановки $\pi$ на множестве
$\{-\infty,0,1,\ldots,6\}$.
\end{theorem}
\proofr Подкод $H_{\alpha^k}\times H_{\alpha^{\pi(k)}}$ для
произвольного $k$ будем называть базовым подкодом кода $D_P$.

Рассмотрим \, $(X,Y)\in H_{\alpha^k}\times H_{\alpha^{\pi(k)}}$ \,
и компоненту \, $R_{\alpha^r,\alpha^s}(D_P,(X,Y))$, \, $r, s\in$
$\in \{-\infty,0,1,\ldots,6\}$. В силу Следствия \ref{corolcom},
три кодовых слова на расстоянии 4 из компоненты
$R_{\alpha^r,\alpha^s}(D_P,(X,Y))$ принадлежат базовому подкоду
$H_{\alpha^k}\times H_{\alpha^{\pi(k)}}$:
$(X+e_{\alpha^r}+e_{\alpha^s}+e_{\alpha^{i_1}}+e_{\alpha^{j_1}},Y),
(X+e_{\alpha^r}+e_{\alpha^s}+e_{\alpha^{i_2}}+e_{\alpha^{j_2}},Y),
(X+e_{\alpha^r}+e_{\alpha^s}+e_{\alpha^{i_3}}+e_{\alpha^{j_3}},Y)$.
В свою очередь, эти векторы находятся на расстоянии 4 от следующих
кодовых слов компоненты $R_{\alpha^r,\alpha^s}(D_P,(X,Y))$:
$$(X+e_{\alpha^{i_1}}+e_{\alpha^{j_1}},Y+e_{\alpha^{i'_1}}+e_{\alpha^{j'_1}})$$
$$(X+e_{\alpha^{i_2}}+e_{\alpha^{j_2}},Y+e_{\alpha^{i'_2}}+e_{\alpha^{j'_2}})$$
$$(X+e_{\alpha^{i_3}}+e_{\alpha^{j_3}},Y+e_{\alpha^{i'_3}}+e_{\alpha^{j'_3}}),$$
для некоторых $i'_1, j'_1, i'_2, j'_2, i'_3, j'_3$. К
рассматриваемым векторам  добавим также кодовое слово
$$(X+e_{\alpha^r}+e_{\alpha^s},Y+e_{\alpha^{r'}}+e_{\alpha^{s'}})$$
кода $D_P$.

В силу Следствия \ref{corolfactor} все эти четыре вектора
принадлежат четырем различным базовым подкодам, отличным от
$H_{\alpha^k}\times H_{\alpha^{\pi(k)}}$, откуда получаем, что
компонента $R_{\alpha^r,\alpha^s}(D_P,(X,Y))$ пересекается с
каждым базовым подкодом.

Теперь рассмотрим произвольный вектор $(X,Y)\in H_{\alpha^k}\times
H_{\alpha^{\pi(k)}}\cap R_{\alpha^r,\alpha^s}(D_P,(X,Y))$ и
покажем, что $X\times H_{\alpha^{\pi(k)}} \subset
R_{\alpha^r,\alpha^s}(D_P,(X,Y))$. Откуда, в силу того, что
компонента всякого кода длины $n=8$ является максимальной и
пересекается с каждым базовым подкодом, получим требуемое.

Пусть $(X,Y)\in H_{\alpha^k}\times H_{\alpha^{\pi(k)}}$,
$(X+e_{\alpha^r}+e_{\alpha^s},Y+e_{\alpha^{r'}}+e_{\alpha^{s'}})\in
H_{\alpha^l}\times H_{\alpha^{\pi(l)}}$. В виду Леммы
\ref{l4graph} для всякого $Y'\in H_{\alpha^{\pi(k)}}$ существует
последовательность $Y, Y^1, Y^2,\ldots, Y^s=Y'$ кодовых слов
$H_{\alpha^{\pi(k)}}$ и $H_{\alpha^{\pi(l)}}$, соседние слова
которой различаются в двух координатных позициях:
$$Y, \,\,
Y^1=Y+e_{\alpha^{a_1}}+e_{\alpha^{b_1}}\in H_{\alpha^{\pi(l)}},$$
$$Y^2=Y^1+e_{\alpha^{a_2}}+e_{\alpha^{b_2}}\in H_{\alpha^{\pi(k)}}$$
$$\ldots$$
$$Y'=Y^{s-1}+e_{\alpha^{a_s}}+e_{\alpha^{b_s}}\in
H_{\alpha^{\pi(k)}}.$$

Всякой такой последовательности отвечает путь в компоненте
$R_{\alpha^r,\alpha^s}(D_P,(X,Y))$ из $(X,Y)$ в $(X,Y')$ в базовых
подкодах $H_{\alpha^l}\times H_{\alpha^{\pi(l)}}$ и
$H_{\alpha^k}\times H_{\alpha^{\pi(k)}}$:

$$(X,Y)\in H_{\alpha^k}\times
H_{\alpha^{\pi(k)}},$$
$$(X+e_{\alpha^r}+e_{\alpha^s},Y^1)=(X+e_{\alpha^r}+e_{\alpha^s},Y+e_{\alpha^{a_1}}+e_{\alpha^{b_1}})\in H_{\alpha^l}\times H_{\alpha^{\pi(l)}}$$
$$(X,Y^2)=(X,Y^1+e_{\alpha^{a_2}}+e_{\alpha^{b_2}})\in H_{\alpha^k}\times H_{\alpha^{\pi(k)}}$$
$$\ldots$$
$$(X,Y')=(X+e_{\alpha^r}+e_{\alpha^s},Y^{s-1}+e_{\alpha^{a_s}}+e_{\alpha^{b_s}})\in H_{\alpha^k}\times H_{\alpha^{\pi(k)}},$$
что доказывает  $X\times H_{\alpha^{\pi(k)}} \subset
R_{\alpha^r,\alpha^s}(D_P,(X,Y))$. 
Максимальность любой ком\-поненты вытекает из того факта, что
каждый код пересекается с базовым подкодом.

\eproof

Заметим, что аналогичный подход (анализ графа расстояний два
 $G_{k,l}$, индуцированного двумя кодами $C_k$ и $C_l$ из разбиения),
может быть применен для исследования структуры компонент кодов
Соловьевой, полученных из других разбиений. Похожими
рассуждениями, например, можно показать, что коды Соловьевой длины
15, полученные из разбиения под номером 1 работы \cite{P2000}
имеют
компоненты мощности хотя бы 512, а для отдельных 
подстановок  и максимальные компоненты по любым однородным
координатам.

\section{Свойства максимальных компонент}

В данном разделе отметим некоторые комбинаторные аспекты
компонент, дающих дополнительную мотивацию для исследования этих
нетривиальных структур. Рассмотрим случай $i$-компонент и
совершенных кодов (длина $n=2^r-1$).
 В этом случае $i$-компонента
определяется как совокупность кодовых слов совершенного кода $D$,
которые можно соединить путем с кодовым словом $X$, каждая
последовательная пара которых состоит из $i$-смежных кодовых слов.

Разбиение вершин графа на \, $m$ \, множеств-цветов,
занумерованных числами $\{1,\ldots,m\}$ назовем {\it совершенной
$m$-раскраской}, если для произвольных $i,j$ всякая вершина цвета
$i$ смежна с $A_{ij}$ вершинами цвета $j$. Квадратная матрица
$A=(A_{ij})$ порядка $m$ называется {\it матрицей параметров}
совершенной раскраски. В случае, когда матрицу параметров
переименованием цветов раскраски можно привести к
трехдиагональному виду, цвет под номером 1 называется {\it
полностью регулярным кодом}. В дальнейшем будем рассматривать лишь
совершенные раскраски двоичных графов Хэмминга. Полностью
регулярные коды были введены Дельсартом \cite{Delsarte} как
обобщения совершенных кодов, сохраняющие их самые замечательные
алгебро-комбинаторные свойства. {\it Радиусом покрытия} $\rho(C)$
кода $C$ назовем максимально возможное расстояние от произвольного
вектора из ${\textbf F}^n$ до кодовых вершин кода $C$.

Кодовое слово совершенного кода называется $i$-{\it четным}
($i$-{\it нечетным}), если носитель кодового слова имеет нечетный
(четный) вес и содержит $i$ (не содержит $i$) или, имея четный
(нечетный) вес, не содержит $i$ (содержит $i$). Прежде всего
покажем связь $i$-четных вершин и максимальных компонент с
совершенными раскрасками и совершенными кодами.

\begin{prop}\label{rhocomp} Пусть $I$ и $I'$ --
множества $i$-четных и $i$-нечетных кодовых слов совершенного кода
$C$ соответственно. Тогда $\rho(I)=\rho(I')$=3 и множество
векторов на расстоянии 3 от множества $I$ есть $I'\cup I'+e_i$.
\end{prop}
\proofr Очевидно,
 что для всякого i-четного кодового слова $x$
найдется $(n-3)(n-1)/6$  кодовых слов на расстоянии 3, являющихся
i-нечетными. Следовательно $d(I,I')=3$.

Рассмотрим произвольный $y$, такой что $d(y,I)\geq 2$, $y \notin
C$. Тогда $y$ покрывается кодовым словом $u\in I'=C\setminus I:
y=u+e_j$. Имеем два случая: $j=i$, тогда $y\in I'+e_i$, но,
поскольку код $I\cup I'+e_i$ -- совершенный с множеством
$i$-четных и $i$-нечетных вершин $I$ и $I'+e_i$,  заключаем, что
$d(y,I)=3$. Если $j\neq i$, тогда найдется кодовое слово $z$,
отличающееся от $x$ в позициях $j, a, b$, причем $a$ и $b$ не
равны $i$, откуда получаем $z\in I$, $d(y,z)=2$. \eproof

Несложно видеть, что радиус покрытия минимальной компоненты
$R=\{(x,|x|,x):x \in F_2^{(n-1)/2}\}$ равен $(n-1)/2$ (в коде
Хэмминга длины $n$ для любого допустимого $n\geq 7$ существует
минимальная компонента, состоящая только из кодовых слов веса
$(n-1)/2$ и $(n+1)/2$). Естественно предположить, что компоненты,
отличные от минимальной и максимальной, имеют радиусы покрытия в
пределах от $4$ до $(n-3)/2$.

 {\bf
Нерешенная проблема.}  Пусть $R$ -- компонента с $\rho=3$.
Является ли она необходимо максимальной?

Пусть $I$, $I'$ есть множества $i$-четных и $i$-нечетных кодовых
слов некоторого двоичного совершенного кода $C$. Рассмотрим
раскраску графа Хэмминга в следующие 6 цветов: $I$, $I+e_i$,
$I_1$, $I'_1$, $I'+e_i$, $I'$, где $I_1$ и $I'_1$ есть вершины на
расстоянии 1 от $I\cup I+e_i$ и $I'\cup I'+e_i$ соответственно.
\begin{corol}
Раскраска графа Хэмминга в цвета $I, I', I_1, I'_1, I+e_i, I'+e_i$
является совершенной с матрицей параметров

\begin{equation}\label{Matrix}
\begin{tabular}{|c|c|c|c|c|c|c|}
  \hline
    & I & $I+e_i$ & $I_1$ & $I'_1$ & $I'+e_i$ & I' \\
 \hline
   I & 0 & 1 & n-1 & 0 & 0 & 0 \\
  $I+e_i$ & 1 & 0 & n-1 & 0 & 0 & 0 \\
  $I_1$ & 1 & 1 & 1 & n-3 & 0 & 0 \\
  $I'_1$ & 0 & 0 & n-3 & 1 & 1 & 1 \\
  $I'+e_i$ & 0 & 0 & 0 & n-1 & 0 & 1 \\
  $I'$ & 0 & 0 & 0 & n-1 & 1 & 0 \\
  \hline
\end{tabular} \,\, .
\end{equation}

Раскраска графа Хэмминга в цвета $I\cup I+e_i, I_1, I'_1, I'\cup
I'+e_i$ является совершенной с матрицей параметров

$$\begin{tabular}{|c|c|c|c|c|}
  \hline
    & $I \cup I+e_i$ & $I_1$ & $I'_1$ & $I'+e_i \cup I'$ \\
 \hline
   $I\cup I+e_i$ & 1 & n-1 & 0 & 0 \\
   $I_1$ & 2 & 1 & n-3 & 0 \\
  $I'_1$ & 0 & n-3 & 1 & 2 \\
  $I'+e_i \cup I'$ & 0 & 0 & n-1 & 1 \\
    \hline
\end{tabular} \,
\, .$$
\end{corol}

{\bf Нерешенная проблема.} Пусть существует совершенная раскраска
графа Хэмминга с матрицей параметров (\ref{Matrix}). Существует ли
такое $i$, что цвет $I$ есть множество $i$-четных вершин
некоторого совершенного кода.

\begin{corol}
Всякая максимальная компонента либо не достраивается максимальной
компонентой до совершенного кода, либо достраивается ровно двумя
способами.
\end{corol}
\proofr Если $I$ и $I'$-множества $i$-четных и $i$-нечетных
кодовых слов кода $C$, причем $I'$ является гигантской
компонентной, то множество $I$ может быть достроено лишь одной
гигантской компонентой до другого совершенного кода, это
компонента $I'+e_i$. Пусть $I''$-гигантская компонента,
дополняющая $I$ до совершенного кода, отличная от $I'$ и $I'+e_i$.
Тогда в силу Утверждения \ref{rhocomp}, $I''\subset I\cup I'$, а
существование пары кодовых слов $x\in I''\cap I'$ и $y\in I''\cap
I'+e_i$, различающихся в трех координатных позициях, включающих
$i$, невозможно.

\eproof

\begin{corol}
Пусть $M$ -- число  попарно неизоморфных совершенных кодов длины
$n$, состоящих из двух максимальных компонент по некоторому
направлению $i$. Тогда найдется хотя бы $M/2$ попарно неизоморфных
максимальных компонент, вложимых в совершенные коды.
\end{corol}

В случае $n=15$ лишь 30 классов изоморфизма из 5983 неизоморфных
классов совершенных кодов не могут быть представлены как
объединение двух максимальных компонент \cite{Pot}. Этот
эмпирический факт подтверждает актуальность исследования таких
объектов как максимальные компоненты совершенных кодов. Как
следствие имеем хотя бы 2977 попарно неизоморфных максимальных
компонент, вложимых в совершенные коды длины n=15.

В заключение авторы выражают благодарность Олли Поттонену за
предоставление информации о структуре компонент совершенных кодов
длины 15.

\end{document}